\documentclass[10pt]{amsart}

\usepackage{amsmath}
\usepackage{amsfonts}
\usepackage{amssymb}
\usepackage{xcolor}
\usepackage{enumerate}
\usepackage{mathrsfs}
\usepackage{lineno}
\usepackage{tikz}
\usepackage{tikz-cd}

\newtheorem{theorem}{Theorem}[section]
\newtheorem{lemma}[theorem]{Lemma}
\newtheorem{proposition}[theorem]{Proposition}
\newtheorem{corollary}[theorem]{Corollary}
\theoremstyle{definition}
\newtheorem{definition}[theorem]{Definition}

\theoremstyle{remark}

\numberwithin{equation}{section}

\newtheorem{question}[theorem]{Question}

\def\cf{\mathop{\rm cf}\nolimits}
\def\cl{\mathop{\rm cl}\nolimits}

\def\dom{\mathop{\rm dom}\nolimits}

\def\Int{\mathop{\rm int}\nolimits}
\def\ran{\mathop{\rm ran}\nolimits}
\def\RO{\mathop{\rm RO}\nolimits}

\let\mathcal\mathscr

\usepackage{hyperref}

\makeatletter
\@namedef{subjclassname@2020}{%
  \textup{2020} Mathematics Subject Classification}
\makeatother

\begin{document}

\title{New results regarding the lattice of uniform topologies on $C(X)$}
\author[R. Pichardo-Mendoza]{Roberto Pichardo-Mendoza}
\address{Departamento de Matem\'aticas, Facultad de Ciencias, Circuito ext. s/n, Ciudad Universitaria, C.P. 04510,  M\'exico, CDMX}
\email{rpm@ciencias.unam.mx}
\urladdr{http://www.matematicas.unam.mx/pmr/}

\author[A. R\'{\i}os-Herrej\'on]{Alejandro R\'{\i}os-Herrej\'on}
\address{Departamento de Matem\'aticas, Facultad de Ciencias, Circuito ext. s/n, Ciudad Universitaria, C.P. 04510,  M\'exico, CDMX}
\email{chanchito@ciencias.unam.mx}
\urladdr{}
\thanks{The research of the second author was supported by CONACyT grant no. 814282.}


\begin{abstract}
  For a Tychonoff space $X$, the lattice $\mathscr U_X$ was introduced in L.A. P\'erez-Morales, G. Delgadillo-Pi\~n\'on, and
  R. Pichardo-Mendoza, \textit{The lattice of uniform topologies on
    $C(X)$}, Questions and Answers in General Topology {\bf 39}
  (2021), 65--71.

  In the present paper we continue the study of $\mathscr U_X$. To be
  specific, the present paper deals, in its first half, with
  structural and categorical properties of $\mathscr U_X$, while in
  its second part focuses on cardinal characteristics of the lattice
  and how these relate to some cardinal functions of the space $X$. 
\end{abstract}

\maketitle


\section{Introduction}
In \cite{perez2021lattice} the authors define, given a completely
regular Hausdorff space $X$, a partially ordered set $(\mathscr U_X , \subseteq)$ (see
Section~\ref{PicRio_secPreli} for details and the corresponding definitions)
which turns out to be a bounded lattice (the {\sl lattice of uniform
  topologies on $C(X)$}). Here we expand some of the results obtained in
that paper and explore new directions. For example, Section~\ref{PicRio_secCatego} is mainly about finding connections
between order-isomorphisms and homeomorphisms, while the last two
sections deal heavily on finding relations between some cardinal
characteristics of $\mathscr U_X$ and highly common cardinal functions
of $X$. 

\section{Preliminaries}\label{PicRio_secPreli}
All topological notions and all set-theoretic notions whose definition
is not included here should be understood as in
\cite{engelking1989general} and \cite{kunen1980set},
respectively. With respect to lattices, we will follow
\cite{larson1975lattice} for notation and results. The same goes for
Boolean algebras and \cite{koppelberg1989handbook1}.

The symbol $\omega$ denotes both, the set of all non-negative integers
and the first infinite cardinal. Also, $\mathbb R$ is the real line
endowed with the Euclidean topology. 

Given a set $S$, $[S]^{<\omega}$ denotes the collection of all finite
subsets of $S$. For a set $A$, the symbol ${}^AS$ is used to represent
the collection of all functions from $A$ to $S$. In particular, for
$f\in{}^AS$, $E\subseteq A$, and $H\subseteq S$ we define
$f[E] := \{f(x) : x\in E\}$ and $f^{-1}[H] := \{x\in A : f(x)\in
H\}$. Moreover, if $y\in S$, $f^{-1}\{y\} := f^{-1}[\{y\}]$.

A nonempty family of sets, $\alpha$, is called {\sl directed} if for
any $A,B\in\alpha$ there is $E\in\alpha$ with $A\cup B\subseteq
E$. For example, $[S]^{<\omega}$ is directed, for any set $S$. 

Assume $X$ is a set. Hence, $\mathcal P(X)$ and $\mathcal D_X$ represent its power set and the collection of all directed subsets of
$\mathcal P(X)$, respectively. In \cite{pichardo2013pseudouniform} the
term {\sl base for an ideal on $X$} was used to refer to members of
$\mathcal D_X$.

Unless otherwise stated, the word {\sl space} means {\sl Hausdorff completely regular space} (i.e., {\sl Tychonoff space}).

Assume $X$ is a space. Then, $\tau_X$ and $\tau_{X}^{*}$ stand,
respectively, for the families of all open and closed subsets of $X$. Moreover, whenever $x\in X$, $\tau_X(x)$ will be the set $\{U\in \tau_X : x\in U\}$. Now, given
$A\subseteq X$, the symbol $\cl_XA$ (or
$\overline A$ when the space $X$ is clear from the context) represents
the closure of $A$ in $X$; similarly, $\mathop{\rm int}_XA$ and
$\mathop{\rm int} A$ will be used to denote the interior of $A$ in
$X$.

$C(X)$ is, as usual, the subset of ${}^X\mathbb R$
consisting of all continuous functions. Now, given $\alpha\in\mathcal D_X$ we generate a
topology on $C(X)$ as follows: a set $U\subseteq C(X)$ is open if and
only if for each $f\in U$ there are $A\in\alpha$ and a real number
$\varepsilon>0$ with \[V(f,A,\varepsilon):=\{g\in C(X) : \forall\,x\in A\ (|f(x) -
g(x)|<\varepsilon)\}\subseteq U.\]The resulting topological space is
denoted by $C_\alpha(X)$. As it is explained in
\cite{pichardo2013pseudouniform}, $C_\alpha(X)$ is a uniformizable
topological space which may not be Hausdorff. In fact, one has the
following result (whose proof can be found in \cite[Proposition~3.1,
p.~559]{pichardo2013pseudouniform}). 

\begin{lemma}\label{caracterizacionDeHausdorff}
  For any space $X$ and $\alpha\in\mathcal D_X$, $C_\alpha(X)$ is
  Hausdorff if and only if $\alpha$ has dense union, i.e., $\overline{\bigcup\alpha}=X$.
\end{lemma}

Given a space $X$, set $\mathscr U_X := \{\tau_{C_\gamma(X)} : \gamma\in\mathscr
  D_X\}$. In order to simplify our writing, for each
  $\alpha\in\mathscr D_X$ we identify the space $C_\alpha(X)$ with its
  topology. Thus, expressions of the form $C_\alpha(X)\in\mathscr U_X$
  will be common in this paper. Also, in those occasions where the ground
  space is clear from the context, we will suppress it from our
  notation, i.e., we will use $C_\alpha$ instead of
  $C_\alpha(X)$. Finally, for any $\alpha,\beta\in\mathscr D_X$, both,
  $C_\alpha(X)\leq C_\beta(X)$ and $C_\alpha\leq C_\beta$, are
  abbreviations of the relation $\tau_{C_\alpha(X)}\subseteq \tau_{C_\beta(X)}$.

It is shown in \cite[Proposition~3.2, p.~67]{perez2021lattice} that
the poset $(\mathscr U_X,\subseteq)$ is a bounded distributive
lattice; to be precise, given $\alpha,\beta\in\mathscr D_X$, the
collections \[\alpha\vee \beta := \{A\cup B : A\in\alpha,\
  B\in\beta\}\qquad\text{and}\qquad \alpha\wedge \beta :=
  \{\overline{A}\cap \overline{B} : A\in\alpha,\ B\in\beta\}\] are
directed and, moreover, $C_{\alpha\vee\beta}$ and
$C_{\alpha\wedge\beta}$ are, respectively, the supremum and infimum of
$\{C_\alpha , C_\beta\}$ in $\mathscr U_X$.

The topologies generated on $C(X)$ by the directed sets
$\{\emptyset\}$, $[X]^{<\omega}$, and $\{X\}$ are denoted by
$C_\emptyset(X)$, $C_p(X)$, and $C_u(X)$, respectively. Let us note that
$C_\emptyset$ is the indiscrete topology on $C(X)$, while $C_p$ and
$C_u$ are the topologies of pointwise and uniform convergence on $C(X)$, respectively.

The result below (see \cite[Theorem~3.4,
p.~560]{pichardo2013pseudouniform} for a proof) will be used several
times in what follows.

\begin{proposition}\label{PerPic_caractMasFina}
  If $X$ is a space and $\alpha,\beta\in\mathcal D_X$, then
  $C_\alpha\leq C_\beta$ if and only if for each $A\in\alpha$ there is $B\in\beta$ with $A\subseteq\overline B$. 
\end{proposition}

We finish this section by mentioning that our notation for topological
cardinal functions follows \cite{hodel1984cardinal}; in particular, all of them are, by
definition, infinite. 

\section{Some structural and categorical results}\label{PicRio_secCatego}
We begin by improving the result presented in \cite[Proposition~3.2, p.~67]{perez2021lattice}.

\begin{proposition}\label{PerPic_retCompleta}
  For any space $X$, $\mathscr U_X$ is a complete lattice.
\end{proposition}

\begin{proof}
  Given an arbitrary set $\mathscr S\subseteq\mathscr D_X$, define $\mathscr A := \{C_\delta : \delta\in\mathscr S\}$.

  By letting $\alpha$ be the family
  of all sets of the form $\bigcup\mathscr E$, where $\mathscr
  E\subseteq\bigcup\mathscr S$ is finite, we obtain $\alpha\in\mathscr
  D_X$. Also, the fact that $\delta\subseteq\alpha$, whenever
  $\delta\in\mathscr S$, implies (see
  Proposition~\ref{PerPic_caractMasFina}) that $C_\alpha$ is an upper
  bound for $\mathscr A$.

  Now, assume that $\gamma\in\mathscr D_X$ is such that $C_\gamma$ is
  an upper bound for $\mathscr A$. In order to show that $C_\alpha\leq C_\gamma$, fix
  $A\in\alpha$. There is a finite set
  $\mathscr E\subseteq\bigcup\mathscr S$ satisfying $A =
  \bigcup\mathscr E$. According to
  Proposition~\ref{PerPic_caractMasFina}, for each $E\in\mathscr E$
  there exists $E^*\in\gamma$ with $E\subseteq\overline{E^*}$. Since
  $\gamma$ is directed, $\bigcup\{E^* : E\in\mathscr E\}\subseteq G$
  for some $G\in\gamma$ and, consequently, $A\subseteq\overline
  G$. In other words, $C_\alpha\leq C_\gamma$. 

  From the previous paragraphs we conclude that any subset of $\mathscr
  U_X$ has a supremum in $\mathscr U_X$. Now, regarding infima, let us
  observe that the infimum of $\emptyset$ in $\mathscr U_X$ is
  $C_u$. Thus, we will suppose that $\mathscr S$ is non-empty.

  Denote by $\mathscr E$ the set of all choice functions of $\mathscr
  S$, i.e., $e\in\mathscr E$ if and only if $e:\mathscr
  S\to\bigcup\mathscr S$ and $e(\delta)\in\delta$, for all
  $\delta\in\mathscr S$. Now, for each $e\in\mathscr E$, set \[\widetilde
  e := \bigcap\{\overline{e(\delta)} : \delta\in\mathscr S\}.\] We
claim that if $\beta := \{\widetilde e : e\in\mathscr E\}$, then
$C_\beta$ is the infimum of $\mathscr A$.

To show that $\beta$ is directed, consider $d,e\in\mathscr E$. Since, for any
$\delta\in\mathscr S$, $\delta$ is directed, we deduce that there is a set
$f(\delta)\in\delta$ with $d(\delta)\cup e(\delta)\subseteq
f(\delta)$. This produces $f$, a choice function of $\mathscr S$, in
such a way that $\widetilde d\cup\widetilde e\subseteq\widetilde f$.

The fact that $C_\beta$ is a lower bound for $\mathscr A$ follows from
the observation that for each $e\in\mathscr E$ and $\delta\in\mathscr
S$, $\widetilde e\subseteq\overline{e(\delta)}$.

Finally, let $\gamma\in\mathscr D_X$ be such that $C_\gamma$ is a
lower bound for $\mathscr A$. Fix $G\in\gamma$. Then, for any
$\delta\in\mathscr S$ there is $e(\delta)\in\delta$ with
$G\subseteq\overline{e(\delta)}$. As a consequence, we obtain $e$, a choice
function of $\mathscr S$, with $G\subseteq\widetilde e$.
\end{proof}

As in \cite{larson1975lattice}, we will use the symbol $\Sigma(E)$ to
the represent the collection of all topologies on a fixed set $E$. It
is well-known that when we order $\Sigma(E)$ by direct inclusion, the
resulting structure is a complete lattice. In
particular, the supremum of $\mathscr A\subseteq\Sigma(E)$ is the
topology on $E$ {\sl generated by} $\bigcup\mathscr A$ (i.e., it has the
collection $\bigcup\mathscr A$ as a subbase).

Clearly, $\mathscr U_X$ is a suborder of $\Sigma(C(X))$. Thus, a
natural question is, given a family $\mathscr A\subseteq\mathscr U_X$,
is the supremum (respectively, infimum) of $\mathscr A$ as calculated in
$\mathscr U_X$ the same as the supremum (respectively, infimum) of
$\mathscr A$ as obtained in $\Sigma(C(X))$? We have a positive answer
for suprema.

\begin{corollary}\label{PicRio_coro:supremoBien}
  If $X$ is a space and $\mathscr A\subseteq\mathscr U_X$, then
  $\bigvee\mathscr A$, the supremum of $\mathscr A$ in
  $\mathscr U_X$, is the topology on $C(X)$ which has
  $\bigcup\mathscr A$ as a subbase. 
\end{corollary}
\begin{proof}
  Fix $\mathscr S\subseteq\mathscr D_X$ in such a way that $\mathscr A
  = \{C_\beta : \beta\in\mathscr S\}$ and denote by $\sigma$ the
  topology on $C(X)$ generated by $\bigcup\mathscr A$. Since
  $\bigvee\mathscr A$ is an upper bound of $\mathscr A$ in $\Sigma(C(X))$, we obtain
  $\sigma\subseteq\bigvee\mathscr A$.

  Now, let $f\in U\in\bigvee\mathscr A$
  be arbitrary. According to the proof of
  Proposition~\ref{PerPic_retCompleta}, there are $\varepsilon>0$ and
  $\mathscr E$, a finite subset of $\bigcup\mathscr S$, with
  $V(f,A,\varepsilon)\subseteq U$, where $A := \bigcup\mathscr
  E$. When $\mathscr E=\emptyset$, we deduce that $U =
  C(X)\in\sigma$. Hence, let us assume that $\mathscr
  E\neq\emptyset$.

  For each $E\in\mathscr E$ let $\beta(E)\in\mathscr S$ be such that
  $E\in\beta(E)$. By setting $\mathscr W := \{\Int_{C_{\beta(E)}} V(f,E,\varepsilon):
  E\in\mathscr E\}$ we produce a finite subset of $\bigcup\mathscr A$ which
  satisfies $f\in\bigcap\mathscr W\subseteq
  V(f,A,\varepsilon)\subseteq U$. In conclusion, $\bigvee\mathscr
  A\subseteq\sigma$. 
\end{proof}

Recall that if $E$ is a set and
$\sigma,\tau\in\Sigma(E)$, the infimum of $\{\sigma,\tau\}$ in
$\Sigma(E)$ is $\sigma\cap\tau$; consequently, for any space $X$ and
$\alpha,\beta\in\mathscr D_X$, $C_\alpha \wedge C_\beta\subseteq
C_\alpha\cap C_\beta$. Now, assume that $X$ is a non-empty space which
is {\sl resolvable} (i.e., it can be written as the union of two
disjoint dense subsets of it). In \cite[Proposition~4.5,
p.~69]{perez2021lattice}, it is shown that there are two Hausdorff
topologies $\sigma,\tau\in\mathscr U_X$ with $\sigma\wedge\tau =
C_\emptyset$. Consequently, $\sigma\cap\tau$ is a $T_1$ topology, but
$\sigma\wedge\tau$ fails to be $T_0$. Hence, the question
posed in the paragraph preceding
Corollary~\ref{PicRio_coro:supremoBien} has a negative answer for
infima.

\begin{question}
  Given a space $X$, find conditions on $\alpha,\beta\in\mathscr{D}_X$
  in order to obtain $C_\alpha \wedge C_\beta = C_\alpha \cap C_\beta$.
  \end{question}

As in \cite{perez2021lattice}, the symbol $\mathscr C_X$ represents the collection of all members of $\mathscr U_X$ which have a
complement in $\mathscr U_X$. Thus, from the fact that $\mathscr U_X$ is a bounded distributive lattice, we deduce that $\mathscr U_X$ is a Boolean algebra if and only if $\mathscr U_X=\mathscr C_X$. Our next result shows that this condition is attained only in trivial cases.

\begin{proposition}
  For any space $X$, $\mathscr U_X$ is a Boolean algebra if and only
  if $X$ is finite. 
\end{proposition}

\begin{proof} Firstly observe that, in virtue of
  \cite[Proposition~3.3, p.~68]{perez2021lattice}, we only need to
  show that $X$ is a finite space if and only if for each $\alpha\in
  \mathscr D_X$ there is $E\in\alpha$ with $\overline{E}\in\tau_X$ and
  $\bigcup\alpha\subseteq\overline{E}$. Now, evidently any finite $X$
  satisfies the latter condition. For the converse let us assume that
  $X$ is infinite. Since $X$ is Hausdorff, there is $\{U_n :
  n<\omega\}$, a family of non-empty open subsets of $X$, with
  $U_m\cap U_n = \emptyset$, whenever $m<n<\omega$. By setting $\alpha :=
  \left\{\bigcup_{k=0}^n U_k : n<\omega\right\}$ we obtain a member of
  $\mathscr D_X$ in such a way that, for each $E\in\alpha$, there is
  $m<\omega$ with $U_m\cap E = \emptyset$ and thus,
  $\bigcup\alpha\not\subseteq\overline{E}$. 
\end{proof}

For our next results we will need some auxiliary concepts. First of
all, assume that $f$ is function from the space $X$ into a space
$Y$. One easily verifies that for any $\alpha\in\mathscr D_X$ the
family \[f^*\alpha := \{f[A] : A\in\alpha\}\] belongs to $\mathscr
D_Y$ and so, we have the following notion (recall that
  for any space $Z$ and $\gamma\in\mathscr D_Z$ we are identifying the
  space $C_\gamma(Z)$ with its topology).

\begin{definition}\label{def_f-induced}
  If $X$, $Y$, and $f$ are as in the previous paragraph, the phrase
  {\sl $\varphi$ is the $f$-induced relation} means that \[\varphi =
    \{(C_\alpha(X),C_{f^*\alpha}(Y)) : \alpha\in\mathscr
    D_X\}\subseteq\mathscr U_X\times\mathscr U_Y.\]
\end{definition}

With the notation used above, the domain of $\varphi$, ${\mathop{\rm
    dom}}(\varphi)$, is equal to $\mathscr U_X$ and its range,
${\mathop{\rm ran}}(\varphi)$, is a subset of $\mathscr U_Y$.

\begin{proposition}
  If $X$ and $Y$ are spaces and $f:X\to Y$, then $f$ is continuous if
  and only if $\varphi$, the $f$-induced relation, is an order-preserving function. 
\end{proposition}
\begin{proof}
  Let us begin by assuming that $f$ is continuous and prove the
  statement below.
  \begin{align}
    \label{PicRio_eq:PropiedadesCategoricas}
    \forall\,\alpha,\beta\in\mathscr D_X\ (C_\alpha\leq C_\beta\ \to\
    C_{f^*\alpha}\leq C_{f^*\beta}).
  \end{align}

  Given $\alpha,\beta\in\mathscr D_X$ with $C_\alpha\leq C_\beta$, fix
  $A\in f^*\alpha$. There is $B\in\alpha$ with $A = f[B]$ and so (see
  Proposition~\ref{PerPic_caractMasFina}), for some $E\in\beta$,
  $B\subseteq{\mathop{\rm cl}}_X E$. Finally, $f$'s continuity
  produces $A = f[B]\subseteq f[{\mathop{\rm cl}}_X
  E]\subseteq{\mathop{\rm cl}}_Y f[E]$ and, clearly, $f[E]\in
  f^*\beta$.

  The final step for this implication is to note that the properties
  required for $\varphi$ are consequences of
  (\ref{PicRio_eq:PropiedadesCategoricas}).

  Suppose that $\varphi$ is an order-preserving function and fix
  $A\subseteq X$. According to Proposition~\ref{PerPic_caractMasFina},
  $C_{{\mathop{\rm cl}}_X A}\leq C_A$ and so, \[C_{f[{\mathop{\rm
          cl}}_X A]} = \varphi(C_{{\mathop{\rm cl}}_X
      A})\leq\varphi(C_A) = C_{f[A]},\] i.e., $f[{\mathop{\rm cl}}_X
  A]\subseteq{\mathop{\rm cl}}_Yf[A]$. 
\end{proof}

For the rest of the paper, given a space $X$, a point $x\in X$, and
a set $A\subseteq X$, we use the symbols $C_x(X)$ and $C_A(X)$ to
represent the topological spaces $C_{\{\{x\}\}}(X)$ and
$C_{\{A\}}(X)$, respectively. As expected, if the space $X$ is clear
from the context, we only write $C_x$ and $C_A$; also, as we have
done before, $C_x$ and $C_A$ are, as well, the topologies of the
corresponding spaces. 

A function $f$ from the space $X$ into the space $Y$ is called
{\sl open onto its range} if, for any $U\in\tau_X$,
$f[U]\in\tau_{f[X]}$. Note that if $f$ is one-to-one, then $f$ is open
onto its range if and only if $f$ is {\sl closed onto its range}
(i.e., whenever $G$ is a closed subset of $X$, $f[G]$ is a closed
subset of the subspace $f[X]$). 

\begin{proposition}
  Assume $X$ and $Y$ are spaces. For any $f : X\to Y$, the following
  are equivalent.
  \begin{enumerate}
  \item $f$ is one-to-one and open onto its range.
  \item $\varphi^{-1}$, the inverse relation of the $f$-induced
    relation, is an order-preserving function. 
  \end{enumerate}
\end{proposition}
\begin{proof}
  Observe that for the implication $(1)\to(2)$, it suffices to prove that the statement
    \begin{align}
    \label{PicRio_eq:PropiedadesCategoricas2}
    \forall\,\alpha,\beta\in\mathscr D_X\ (C_{f^*\alpha}\leq C_{f^*\beta} \ \to\
C_\alpha\leq C_\beta )
    \end{align}
    follows from $(1)$. Thus, suppose $(1)$ and fix $\alpha,\beta\in\mathscr D_X$ with $C_{f^*\alpha}\leq C_{f^*\beta}$. Given
    $A\in\alpha$, Proposition~\ref{PerPic_caractMasFina} guarantees
    the existence of $B\in\beta$ with $f[A]\subseteq{\mathop{\rm
        cl}}_Y f[B]$, i.e., $A\subseteq f^{-1}[{\mathop{\rm cl}}_Y
    f[B]]$. Thus, we only need to show that $f^{-1}[{\mathop{\rm
        cl}}_Y f[B]]\subseteq {\mathop{\rm cl}}_X B$. If $x\in
    f^{-1}[{\mathop{\rm cl}}_Y f[B]]$ and $U\in\tau_X(x)$ are
    arbitrary, then $f(x)\in f[X]\cap{\mathop{\rm cl}}_Y f[B] =
    {\mathop{\rm cl}}_{f[X]} f[B]$ and $f[U]\in\tau_{f[X]}(f(x))$;
    consequently, $f[U]\cap f[B]\neq\emptyset$. Since $f$ is
    one-to-one, $f[U\cap B]\neq\emptyset$ and so, $U\cap
    B\neq\emptyset$, as required.

    For the rest of the argument, assume $(2)$. In order to verify
    that $f$ is one-to-one, let $x,y\in X$ be such that $f(x) =
    f(y)$. Hence, $C_{f(x)} = C_{f(y)}$ and, as a consequence, $C_x =
      \varphi^{-1}(C_{f(x)}) = \varphi^{-1}(C_{f(y)}) = C_y$. The use
      of Proposition~\ref{PerPic_caractMasFina} produces $x=y$. 

      Given that $f$ is
      one-to-one, we only need to argue that $f$ is closed onto its
      range. Suppose $G$ is a closed subset of $X$. By letting $E := {\mathop{\rm
          cl}}_Y f[G]$ and $A := f^{-1}[E]$, we deduce that $f[A] =
      E\cap f[X] = {\mathop{\rm cl}}_{f[X]} f[G]$. Therefore, $C_{f[A]} \leq C_E \leq C_{f[G]}$ and
      so, $C_A = \varphi^{-1}(C_{f[A]}) \leq \varphi^{-1}(C_{f[G]}) =
      C_G$. Hence, $A\subseteq {\mathop{\rm cl}}_X G = G$ and,
      consequently, ${\mathop{\rm cl}}_{f[X]} f[G] = f[A] \subseteq
      f[G]$, i.e., $f[G]$ is a closed subset of $f[X]$. 
\end{proof}

\begin{proposition}
  If $X$ and $Y$ are spaces and $f:X\to Y$, then $f$ is onto if
  and only if ${\mathop{\rm ran}} (\varphi) = \mathscr U_Y$, where
  $\varphi$ is the $f$-induced relation.
\end{proposition}
\begin{proof}
  When $f$ is onto and $\alpha\in\mathscr D_Y$, the collection $\beta
  := \{f^{-1}[A] : A\in\alpha\}$ belongs to $\mathscr D_X$ and
  $f^*\beta = \alpha$. Thus, $(C_\beta,C_\alpha)\in\varphi$ and so,
  $C_\alpha\in {\mathop{\rm ran}} (\varphi)$. 

  For the remaining implication, fix $y\in Y$ and note that
  $C_y\in\mathscr U_Y = {\mathop{\rm ran}} (\varphi)$, i.e., for some
  $\alpha\in\mathscr D_X$, $(C_\alpha,C_y)\in\varphi$. Now, our
  definition of $\varphi$ produces $\beta\in\mathscr D_X$ with
  $C_\alpha = C_\beta$ and $C_y = C_{f^*\beta}$. Since $C_y\leq
  C_{f^*\beta}$, there is $B\in\beta$ in such a way that $y\in
  {\mathop{\rm cl}}_X f[B]$ and so, $B\neq\emptyset$. From the
  relation $C_{f^*\beta}\leq C_y$ we obtain $f[B]\subseteq
  {\mathop{\rm cl}}_Y\{y\} = \{y\}$ and therefore, $\emptyset\neq
  B\subseteq f^{-1}\{y\}$. 
\end{proof}

Since any topological embedding is a
continuous one-to-one function that is open onto its range, we obtain
the following result.

\begin{corollary}\label{order_embedding}
  If $Y$ is a space which can be embedded into a space $X$, then
  there is an order-embedding from $\mathscr U_Y$ into $\mathscr
  U_X$. In particular, $|\mathscr U_Y|\leq |\mathscr U_X|$.
\end{corollary} 

Assume $X$ and $Y$ are spaces for which there is $\varphi : \mathscr
U_X\to\mathscr U_Y$, an (order) isomorphism. According to
\cite[Proposition~5.1, p.~70]{perez2021lattice}, for each $x\in X$,
$C_x(X)$ is an atom of $\mathscr U_X$ (i.e., a minimal element of
$\mathscr U_X\setminus\{C_\emptyset\}$) and so, $\varphi(C_x(X))$ happens to be an atom of
$\mathscr U_Y$; consequently (see \cite[Proposition~5.1,
p.~70]{perez2021lattice}), there exists a point $y\in Y$ with
$\varphi(C_x(X)) = C_{y}(Y)$. Moreover, as one easily deduces from
Proposition~\ref{PerPic_caractMasFina}, $y$ is the only member of $Y$
with this property. 

\begin{definition}\label{defi_PicRio_induced}
  Let $X$ and $Y$ be a pair of spaces. If $\varphi : \mathscr
  U_X\to\mathscr U_Y$ is an isomorphism, we will say that
  $f:X\to Y$ is the  {\sl $\varphi$-induced function} if
  \begin{eqnarray}
    \label{PicRio_eq:1}
    \text{for each }x\in X,\ \varphi(C_x(X)) = C_{f(x)}(Y).
  \end{eqnarray}
\end{definition}

Observe that if $f$ is a homeomorphism from a space $X$ onto a space
$Y$ and $\varphi$ is the $f$-induced relation, the previous
results imply that $\varphi$ is an isomorphism. Now, when $g$ is the
$\varphi$-induced function, we obtain that, for each $x\in
X$, \[\varphi(C_x) = C_{f^*\{\{x\}\}} = C_{f(x)}\qquad\text{and}\qquad
  \varphi(C_x) = C_{g(x)},\] i.e., $f(x) = g(x)$. In conclusion,
$f=g$. Hence, the following is a natural question. 

\begin{question}\label{PicRio_pregElOtroCamino}
  Assume $X$ and $Y$ are spaces for which there is an isomorphism $\varphi
  : \mathscr U_X\to\mathscr U_Y$. If $f$ is the $\varphi$-induced
  function and $\psi$ is the $f$-induced relation, do we get $\varphi
  = \psi$?
\end{question}

With the idea in mind of giving a positive answer to this question for
a class of spaces (zero-dimensional spaces), we will present some auxiliary results.

\begin{lemma}\label{PicRio_fEsBiyectiva}
  Assume $\varphi : \mathscr U_X\to\mathscr U_Y$ is an isomorphism,
  where $X$ and $Y$ are spaces. If $f$ is the $\varphi$-induced
  function, then the following statements hold.
  \begin{enumerate}
  \item $f$ is a bijection and $f^{-1}$ is the $\varphi^{-1}$-induced
    function.
  \item If $A\subseteq X$ and $\beta\in\mathscr D_Y$ satisfy
    $\varphi(C_A(X)) = C_\beta(Y)$, then $f[\cl_X A]\subseteq\bigcup\overline{\beta}$. 
  \end{enumerate}
\end{lemma}
\begin{proof}
  For (1), let $g$ be the $\varphi^{-1}$-induced function. Given $x\in
  X$, the relation $\varphi(C_x) = C_{f(x)}$ implies that $C_x =
  \varphi^{-1}(C_{f(x)}) = C_{g(f(x))}$ and so, $g\circ f$ is the
  identity function on $X$. Similarly, $f\circ g$ is the identity
  function on $Y$.

  Given $x\in\overline A$, Proposition~\ref{PerPic_caractMasFina}
  produces $C_x\leq C_A$ and so, $C_{f(x)} = \varphi(C_x)\leq\varphi(C_A) = C_\beta$;
  hence, $f(x)\in\bigcup\overline\beta$. 
\end{proof}

\begin{proposition}\label{PicRio_equivHomeoIsom}
  Let $X$ and $Y$ be spaces in such a way that there is an isomorphism
  $\varphi:\mathscr U_X\to\mathscr U_Y$. Denote by $f$ the
  $\varphi$-induced function and consider the following statements.
  \begin{enumerate}
  \item $\varphi$ is the $f$-induced relation.
  \item For any $A\subseteq X$, $\varphi(C_A(X)) = C_{f[A]}(Y)$.
  \item Whenever $G$ is a closed subset of $X$, $\varphi(C_G(X)) = C_{f[G]}(Y)$.
  \end{enumerate}
  Then, (1) is equivalent to (2) and if $f$ is continuous, (2) and (3)
  are equivalent. 
\end{proposition}
\begin{proof}
  The implications (1)$\to$(2) and (2)$\to$(3) are immediate. On the
  other hand, it follows from the work done in the first paragraphs of
  the proof of Proposition~\ref{PerPic_retCompleta} that, for any
  $\alpha\in\mathscr D_X$, \[C_\alpha = \bigvee\{C_A :
    A\in\alpha\}\qquad\text{and}\qquad C_{f^*\alpha} =
    \bigvee\{C_{f[A]} : A\in\alpha\};\] therefore, by assuming (2) we
  obtain \[\varphi(C_\alpha) = \bigvee\{\varphi(C_A) : A\in\alpha\} =
    \bigvee\{C_{f[A]} : A\in\alpha\} = C_{f^*\alpha},\] i.e., (1)
  holds. 

  Now suppose $f$ is continuous and (3) is true. In order to prove (2),
  fix $A\subseteq X$ and set $G := \overline A$. According to
  Proposition~\ref{PerPic_caractMasFina}, $C_A =
  C_G$ and, consequently, $\varphi(C_A) = \varphi(C_G) = C_{f[G]}$. From
  the relation $f[A]\subseteq f[G]$ we deduce that $C_{f[A]}\leq
  C_{f[G]}$. The continuity of $f$ produces $f[G]\subseteq
  \overline{f[A]}$ and so, $C_{f[G]}\leq C_{f[A]}$. In conclusion,
  $\varphi(C_A) = C_{f[G]} = C_{f[A]}$, as needed. 
\end{proof}

Recall that for any space $Z$, $\mathop{\rm CO}(Z)$ is the collection
of all subsets of $Z$ which are closed and open in $Z$. Consequently,
$Z$ is {\sl zero-dimensional} when $\mathop{\rm CO}(Z)$ is a base for
$Z$. 

\begin{lemma}\label{PicRio_todoBienConClopens}
  Assume $X$ and $Y$ are spaces for which there is $\varphi :\mathscr
  U_X\to\mathscr U_Y$, an isomorphism. If $f$ is the $\varphi$-induced
  function, the following statements hold. 
  \begin{enumerate}
  \item For each $A\in\mathop{\rm CO}(X)$, $f[A]\in\mathop{\rm CO}(Y)$ and
    $\varphi(C_A(X)) = C_{f[A]}(Y)$.
  \item If $Y$ is zero-dimensional, $f$ is continuous. 
  \end{enumerate}
\end{lemma}
\begin{proof}
  Given $A\in\mathop{\rm CO}(X)$, the proof of \cite[Proposition~3.3,
  p.~68]{perez2021lattice} shows that $C_A$ and $C_{X\setminus A}$ are
complements of each other in $\mathscr U_X$ and so, $\varphi(C_A)$ and
$\varphi(C_{X\setminus A})$ have the same relation in $\mathscr
U_Y$. Then, according to \cite[Proposition~5.3,
p.~70]{perez2021lattice}, there exists $B\in\mathop{\rm CO}(Y)$ with
$\varphi(C_A) = C_B$ and $\varphi(C_{X\setminus A}) = C_{Y\setminus
  B}$. From Lemma~\ref{PicRio_fEsBiyectiva}(2), $f[\overline
A]\subseteq \overline B$ and $f[\overline{X\setminus A}] \subseteq
\overline{Y\setminus B}$, i.e., $f[A]\subseteq B$ and $Y\setminus
B\supseteq f[X\setminus A] = Y\setminus f[A]$. Thus, $f[A] =
B$.

For the second part, fix $B\in\mathop{\rm
  CO}(Y)$. According to Lemma~\ref{PicRio_fEsBiyectiva}(1), $f^{-1}$
is the $\varphi^{-1}$-induced function and so, we can apply part (1)
of this lemma to $f^{-1}$ in order to get $f^{-1}[B]\in\tau_X$. Thus,
the assumption that $\mathop{\rm CO}(Y)$ is a base for $Y$ gives $f$'s continuity. 
\end{proof}

\begin{lemma}\label{PicRio_elLemaQueSirvePara0Dim}
  Let $X$ and $Y$ be spaces, with $X$ zero-dimensional. If $\varphi$
  is an isomorphism from $\mathscr U_X$ onto $\mathscr U_Y$ and $f$
  is the $\varphi$-induced function, then $\varphi(C_G)\leq C_{f[G]}$,
  whenever $G$ is a closed subset of $X$.
\end{lemma}
\begin{proof}
  Given $G$, a closed subset of $X$, there are $\mathscr
  A\subseteq\mathop{\rm CO}(X)$ and $\beta\in\mathscr D_X$ in such a
  way that $G = \bigcap\mathscr A$ and $\varphi(C_G) =
  C_\beta$. Let us argue that
  \begin{align}
    \label{PicRio_ecuDom0Dim}
  \text{for all }A\in\mathscr A\text{ and }B\in\beta,\ B\subseteq f[A].
  \end{align}

  Suppose $A\in\mathscr A$ and $B\in\beta$ are arbitrary. Since $G\subseteq
  A$, we deduce that $C_G\leq C_A$ and, consequently, the use of
  Lemma~\ref{PicRio_todoBienConClopens}(1) gives \[C_\beta = \varphi(C_G)
    \leq \varphi(C_A) = C_{f[A]};\] in particular,
  $B\subseteq\overline{f[A]}$. To complete this part, invoke
  lemmas~\ref{PicRio_fEsBiyectiva}(1) and
  \ref{PicRio_todoBienConClopens}(2) in order to get the continuity of
  $f^{-1}$, i.e., the closedness of $f$.

  From (\ref{PicRio_ecuDom0Dim}) and the fact that $f$ is one-to-one,
  we obtain that, for any $B\in\beta$, \[B\subseteq\bigcap\{f[A] :
    A\in\mathscr A\} = f\left[\bigcap\mathscr A\right] = f[G].\] In
  other words, $C_\beta\leq C_{f[G]}$, as claimed. 
\end{proof}

\begin{proposition}
  Let $X$, $Y$, $\varphi$, $f$, and $\psi$ be as in
  Question~\ref{PicRio_pregElOtroCamino}. If $X$ and $Y$ are
  zero-dimensional, then $\varphi=\psi$. 
\end{proposition}
\begin{proof}
  First of all, lemmas~\ref{PicRio_todoBienConClopens}(2) and
  \ref{PicRio_fEsBiyectiva}(1) guarantee that $f$ is a
  homeomorphism.

  With the idea in mind of verifying condition (3) of 
  Proposition~\ref{PicRio_equivHomeoIsom}, fix $G$, a closed subset of
  $X$. According to Lemma~\ref{PicRio_elLemaQueSirvePara0Dim},
  $\varphi(C_G)\leq C_{f[G]}$. On the other hand, $f[G]$ is a closed
  subset of $Y$ and so, by applying
  Lemma~\ref{PicRio_elLemaQueSirvePara0Dim} to $\varphi^{-1}$ and
  $f^{-1}$, we obtain $\varphi^{-1}(C_{f[G]})\leq C_{f^{-1}[f[G]]} =
  C_G$, i.e., $C_{f[G]}\leq \varphi(C_G)$. Thus, $\varphi(C_G) =
  C_{f[G]}$.

  We conclude that $\varphi$ is the $f$-induced relation or, in other
  words, $\varphi = \psi$. 
\end{proof}

\begin{corollary}\label{PicRio_elCoroCon0Dim}
  Let $X$ and $Y$ be a pair of zero-dimensional spaces. For any
  function $\varphi : \mathscr U_X\to\mathscr U_Y$, the following
  statements are equivalent.
  \begin{enumerate}
  \item $\varphi$ is an isomorphism.
  \item For some homeomorphism $f : X\to Y$, $\varphi$ is the $f$-induced relation.
  \end{enumerate}
\end{corollary}

\begin{question}
  Is the assumption of zero-dimensionality necessary in
  Corollary~\ref{PicRio_elCoroCon0Dim}? To be more precise, are
  there non-homeomorphic spaces $X$ and $Y$ for which the
  lattices $\mathscr U_X$ and $\mathscr U_Y$ are isomorphic?
\end{question}

\section{Some cardinal characteristics}

\begin{definition}
  For a space $X$, set $\mathscr U_X^+ := \mathscr
  U_X\setminus\{C_\emptyset\}$. Also, given a family $\mathscr
  S\subseteq\mathscr U_X^+$, we say that
  \begin{enumerate}
  \item $\mathscr S$ is {\sl an antichain in $\mathscr U_X$} if for any
    $\sigma,\tau\in\mathscr S$, the condition $\sigma\neq\tau$ implies
    that $\sigma\wedge\tau = C_\emptyset$;
  \item $\mathscr S$ is {\sl dense in $\mathscr U_X$} if for each
    $\sigma\in\mathscr U_X^+$ there is $\tau\in\mathscr S$ with
    $\tau\leq\sigma$. 
  \end{enumerate}
\end{definition}

For a space $X$, {\sl the cellularity of $\mathscr U_X$},
$c(\mathscr U_X)$, is the supremum of all cardinals of the form
$|\mathscr W|$, where $\mathscr W$ is an antichain in $\mathscr
U_X$. The {\sl density of $\mathscr U_X$}, $\pi(\mathscr U_X)$, is the
minimum size of a dense subset of $\mathscr U_X$.

\begin{proposition}
  If $X$ is a space, then $c(\mathscr U_X) = \pi(\mathscr U_X) =
  |X|$. 
\end{proposition}
\begin{proof}
  As one easily verifies, $\mathscr A := \{C_{x} : x\in X\}$ is an
  antichain in $\mathscr U_X$. Thus, $|X|\leq c(\mathscr U_X)$. On the
  other hand, if $\alpha\in\mathscr D_X$ satisfies
  $C_\alpha\in\mathscr U_X^+$, then $C_\alpha\not\leq C_\emptyset$,
  i.e., there are $A\in\alpha$ and $z\in A$. Therefore,
  $C_{z}\leq C_\alpha$ and, consequently, $\mathscr A$ is a
  dense subset of $\mathscr U_X$. Hence, $\pi(\mathscr U_X)\leq|X|$.

  In order to prove that $c(\mathscr U_X)\leq\pi(\mathscr U_X)$, let us fix $\mathscr W$, an antichain in $\mathscr
  U_X$, and $\mathscr S$, a dense subset of $\mathscr U_X$. Then,
  there is $e:\mathscr W\to\mathscr S$ such that
  $e(\tau)\leq\tau$, whenever $\tau\in\mathscr W$. Given
  $\sigma,\tau\in\mathscr W$ with $\sigma\neq\tau$, one gets
  $e(\sigma)\wedge e(\tau)\leq\sigma\wedge\tau = C_\emptyset$ and so,
  $e(\sigma)\neq e(\tau)$; in other words, $e$ is one-to-one and, as a
  consequence, $|\mathscr W|\leq|\mathscr S|$. 
\end{proof}

Now we turn our attention to $|\mathscr U_X|$ and $|\mathscr D_X|$,
for an arbitrary space $X$. With this in mind, given a cardinal $\kappa$, let us recursively define $\beth_0(\kappa)
:= \kappa$ and, for each integer $n$, $\beth_{n+1}(\kappa) :=
2^{\beth_n(\kappa)}$.

\begin{proposition}\label{alelemma}
  The following statements hold for any finite space $X$. 
  \begin{enumerate}
  \item When $|X| = 1$,  $|\Sigma(X)|<2^{|X|} < |\mathscr D_X| = \beth_2(|X|)$.
  \item If $X$ has at least two points, then $2^{|X|} \leq |\Sigma(X)|<|\mathscr D_X| < \beth_2(|X|)$.
  \item $|\mathscr U_X|=2^{|X|}$. 
  \end{enumerate}
  \end{proposition}
  \begin{proof}
    If $X$ has exactly one element,
    then \[\Sigma(X)=\{\{\emptyset,X\}\} \quad \text{and}
      \quad\mathscr D_X=\{\emptyset, \{\emptyset\}, \{X\} ,
      \{\emptyset,X\}\}.\]
    
    With respect to (2), since the function $\eta : \mathscr
    P(X)\setminus\{\emptyset\}\to\Sigma(X)$ given by $\eta(A) :=
    \{\emptyset , A , X\}$ is one-to-one, we deduce that $2^{|X|}-1 =
    |\ran(\eta)|\leq|\Sigma(X)|$. Let us fix $p,q\in X$ with $p\neq
    q$. From the fact that $\{\emptyset,\{p\},\{q\},\{p,q\},X\}$ is a member
    of $\Sigma(X)\setminus\ran(\eta)$, it follows that
    $2^{|X|}\leq|\Sigma(X)|$. 

    The relations $\Sigma(X)\subseteq \mathscr D_X$ and $\{X\} \in \mathscr
    D_X \setminus \Sigma(X)$ clearly imply that $|\Sigma(X)|<|\mathscr
    D_X|$. Lastly, the inequality $|\mathscr D_X| < \beth_2(|X|)$
    follows from the facts $\mathscr D_X\subseteq\mathscr P(\mathscr P(X))$ and
    $C_p\vee C_q\in \mathscr P(\mathscr P(X))\setminus \mathscr D_X$.

    In order to prove (3), start by noticing that from $|X|<\omega$
    one gets $C_p=C_u$. Thus, \cite[Proposition~5.2,
    p.~70]{perez2021lattice} implies that $\mathscr P(X)$, ordered by
    direct inclusion, and the closed interval $[C_\emptyset,C_u]$,
    equipped with the order it inherits from $\mathscr U_X$, are
    order-isomorphic. Finally, (1) in \cite[Proposition~3.2,
    p.~67]{perez2021lattice} guarantees that $\mathscr U_X = [C_\emptyset,C_u]$.
  \end{proof}

  Given a space $X$, let us denote by $\RO(X)$ the collection of all regular open subsets of $X$. According to \cite[Theorem~1.37,
p.~26]{koppelberg1989handbook1}, when we order $\mathop{\rm
  RO}(X)$ by direct inclusion, the resulting structure is a complete
Boolean algebra.
  
\begin{proposition}\label{roblemma}
  The following relations hold for any infinite topological space $X$.
  \begin{enumerate}
  \item $|\mathscr D_X| = \beth_2(|X|)$.
  \item $\max\left\{2^{|X|},2^{|\mathop{\rm RO}(X)|}\right\}\leq
    |\mathscr U_X|\leq 2^{o(X)}$, where $o(X) := |\tau_X|$.
  \end{enumerate}
\end{proposition}
\begin{proof}
  The inequality $|\mathscr D_X| \leq \beth_2(|X|)$ follows from the
  relation $\mathscr D_X\subseteq\mathscr P(\mathscr P(X))$. On the
  other hand, according to \cite[Theorem~7.6, p.~75]{jech2003set},
  there are $\beth_2(|X|)$ filters on the set $X$ and, naturally,
  each one of them is a member of $\mathscr D_X$. This proves (1).
  
  With respect to (2), recall that $\tau_{X}^{*}$ is the collection of all closed
  subsets of $X$. Clearly, $|\tau_{X}^{*}|=o(X)$. An
  immediate consequence of Proposition~\ref{PerPic_caractMasFina}
  is that for each $\alpha\in\mathcal D_X$ the family
  $\overline{\alpha} := \{\overline A : A\in\alpha\}$ is a directed
  set and $C_\alpha = C_{\overline{\alpha}}$. Therefore, $\mathscr
  U_X$ is equal to $\{C_\beta : \beta\in\mathscr D_X\ \wedge\
  \beta\subseteq\tau_{X}^{*}\}$, which, in turn, implies that
  $|\mathscr U_X|\leq|\mathscr P(\tau_{X}^{*})|=2^{o(X)}$.
  
  Now, \cite[Proposition~5.2, p.~70]{perez2021lattice} guarantees the
  existence of a one-to-one map from $\mathscr P(X)$ into $\mathscr
  U_X$ and so, $2^{|X|}\leq |\mathscr U_X|$. 
  
  For the remaining inequality we need some notation. First, given a
  finite function $p\subseteq\mathop{\rm RO}(X)\times 2$, set \[p^\sim := p^{-1}\{0\}\cup\{-x
    : x\in p^{-1}\{1\}\},\] where $-x$ is the Boolean complement of $x\in
  \mathop{\rm RO}(X)$. Hence, a set $\mathscr A\subseteq\mathop{\rm
    RO}(X)$ is called {\sl independent} if for any finite function
  $p\subseteq\mathscr A\times 2$ one has $\bigwedge p^\sim\neq\emptyset$. 
  
  The fact that $X$ is an infinite Tychonoff space
  implies that $\mathop{\rm RO}(X)$ is infinite as well and so, by
  Balcar-Fran\v ek's Theorem (see \cite[Theorem~13.6,
  p.~196]{koppelberg1989handbook1}), there is an independent set $\mathscr A\subseteq
  \mathop{\rm RO}(X)$ with $|\mathscr A| = |\mathop{\rm
    RO}(X)|$. 
  
  Let us argue that, for each $d: \mathscr A\to 2$, the
  collection \[\alpha(d) := \left\{\bigvee p^\sim :
      p\in[d]^{<\omega}\right\}\] is a member of $\mathscr D_X$. Indeed,
  if $p,q\in[d]^{<\omega}$, then $r:=p\cup q$ is a finite subset of $d$
  with $\bigvee r^\sim = \left(\bigvee p^\sim\right)\vee\bigvee
  q^\sim$ and since $\mathop{\rm RO}(X)$ is ordered by direct inclusion,
  we conclude that $\bigvee r^\sim$ is an element of $\alpha(d)$ which
  is a superset of $\bigvee p^\sim$ and $\bigvee q^\sim$. \medskip
  
  \noindent{\bf Claim.} If $d,e\in{}^{\mathscr A}2$ and $U\in\mathscr
  A$ satisfy $d(U) =0$ and $e(U)= 1$, then, for any $V\in\alpha(e)$,
  $U\not\subseteq V$. \medskip
  
  Before we present the proof of our Claim, let's assume it holds and
  fix $d,e\in{}^{\mathscr A}2$ with $d\neq e$. Without loss of
  generality, we may assume that, for some $U\in\mathscr A$, $d(U)=0$
  and $e(U)=1$. Thus, $U\in\alpha(d)$ and if $V$ were a member of $\alpha(e)$ with
  $U\subseteq\overline V$, we would get $U = \mathop{\rm int} U
  \subseteq\mathop{\rm int}\overline V = V$, a contradiction to the
  Claim. As a consequence of this argument, we obtain that the
  function from ${}^{\mathscr A}2$ into $\mathscr U_X$ given by
  $d\mapsto C_{\alpha(d)}$ is one-to-one and so, $2^{|\mathop{\rm RO}(X)|} = 2^{|\mathscr A|}\leq |\mathscr U_X|$. 
  
  Suppose $d$, $e$, and $U$ are as in the Claim. Seeking a
  contradiction, let us assume that $U\subseteq\bigvee p^\sim$, for some
  $p\in[e]^{<\omega}$. We affirm that if $q :=
  p\restriction(\dom(p)\setminus\{U\})$ (the restriction of the
  function $p$ to the given set), then
  \begin{align}
    \label{PicRio_eq:1}
    U\subseteq\bigvee q^\sim.
  \end{align}
  Indeed, when $U\notin\dom(p)$, $p=q$. On the other hand, if
  $U\in\dom(p)$, the relation $p\subseteq e$ gives $p(U) = 1$ and so,
  $\bigvee p^\sim = (-U)\vee\bigvee q^\sim$ which, clearly, implies (\ref{PicRio_eq:1}).
  
  Let us define $r : \dom(q)\cup\{U\}\to2$ by $r(V) = 1-q(V)$, whenever
  $V\in\dom(q)$, and $r(U) = 0$. Obviously, $r\subseteq\mathscr A\times
  2$ is a finite function and thus, the independence of $\mathscr A$
  and the De Morgan's laws produce \[\emptyset\neq\bigwedge
    r^\sim=U\wedge\left(-\bigvee q^\sim\right),\] a contradiction to
  (\ref{PicRio_eq:1}). 
\end{proof}

  Let us recall that a {\sl $T_6$-space} (equivalently,  {\sl perfectly
  normal space}) is a Hausdorff normal space in which all open sets
are of type $F_\sigma$.

  \begin{corollary}\label{cardinality_perfectly_normal_U_X}
  If $X$ is a $T_6$-space, then $|\mathscr U_X|=2^{o(X)}$.

\end{corollary}
\begin{proof}
  We only need to mention that, according to \cite[Theorem~10.5,
  p.~40]{hodel1984cardinal}, $|\mathop{\rm RO}(X)|=o(X)$.
\end{proof}

Our next result is a direct consequence of
corollaries~\ref{cardinality_perfectly_normal_U_X} and
\ref{order_embedding} (recall that any infinite Tychonoff space
contains a copy of the discrete space of size $\omega$).

\begin{corollary}\label{cardinality_discrete_U_X}
  If $Y$ is an infinite discrete subspace of a space $X$, $\beth_2(|Y|)
  \leq |\mathscr U_X|$. In particular, when $X$ is infinite,
  $2^{\mathfrak c}\leq |\mathscr U_X|$. 
\end{corollary}

  Standard arguments show that if $X$ is an arbitrary space and $D$ is a dense
  subspace of it, then the function from
  $\mathop{\rm RO}(X)$ into $\mathscr P(D)$ given by $U\mapsto U\cap D$
  is one-to-one. Therefore (recall that $d(X)$ is the density of $X$),
  \begin{align}
    \label{PicRio3_eqROyd}
    \text{for any space }X,\ |\mathop{\rm RO}(X)|\leq2^{d(X)}.
  \end{align}

  Regarding the accuracy of the bounds presented in
  Proposition~\ref{roblemma}(2), we have the result below. 

\begin{proposition}
  The following statements are true.
  \begin{enumerate}
  \item If $X$ is the Moore-Niemytzki plane (see
    \cite[Example~1.2.4, p.~21]{engelking1989general}), then $|X| = |\mathop{\rm RO}(X)| =
    \mathfrak c$ and $o(X) = 2^{\mathfrak c}$.
  \item When $X$ is the Stone-\v Cech compactification of the
    integers, $|\mathop{\rm RO}(X)| = \mathfrak c$ and $|X| = o(X) =
    2^{\mathfrak c}$.
  \item If $X$ is the Arens-Fort space, \cite[Example~1.6.19,
    p.~54]{engelking1989general}, then $|X| = \omega$ and
    $|\mathop{\rm RO}(X)| = o(X) = \mathfrak c$. 
  \end{enumerate}
\end{proposition}
\begin{proof}
  Let us prove (1). Clearly, $|X| = \mathfrak c$. The equality $|\mathop{\rm RO}(X)| = \mathfrak c$
  follows from the facts, (i) property (\ref{PicRio3_eqROyd}) (recall
  that $X$ is separable) and (ii) the canonical base for $X$ consists of
  $\mathfrak c$ many regular open sets.  Note that from (ii) we also
  deduce the relation $o(X)\leq 2^{\mathfrak c}$. Finally, since
  $X\setminus(\mathbb R\times\{0\})$ is
  an open subset of $X$ which is homeomorphic to an open subspace of the
  Euclidean plane, we conclude that $2^{\mathfrak c}\leq o(X)$.

  Suppose $X$ is as in (2). From \cite[Corollary~3.6.12,
  p.~175]{engelking1989general}, $|X| = 2^{\mathfrak c}$. On the other
  hand, the relation $|\mathop{\rm RO}(X)| = \mathfrak c$ is a
  consequence of (\ref{PicRio3_eqROyd}) and the fact that, according
  to Theorem~3.6.13 and Corollary~3.6.12 of
  \cite[p.~175]{engelking1989general}, $X$ is a space of weight
  $\mathfrak c$ possessing a base of closed-and-open sets. This last statement also implies that $o(X)\leq2^{\mathfrak
    c}$. Now, \cite[Example~3.6.18, p.~175]{engelking1989general}
  guarantees that $X$ has a pairwise disjoint family consisting of
  $\mathfrak c$ many non-empty open sets and so, $2^{\mathfrak
    c}\leq o(X)$.

  Finally, when $X$ is as in (3), one clearly gets $|X|=\omega$ and,
  therefore, $o(X)\leq\mathfrak c$. On the other hand, by
  definition, $X$ has a base consisting of $\mathfrak c$ many
  closed-and-open sets; hence, $\mathfrak c \leq |\mathop{\rm
    RO}(X)|\leq o(X)$. 
\end{proof}

In the next section we focus on the problem of calculating $|\mathscr
U_X|$, for some spaces $X$.

\section{The size of $\mathscr U_X$}
Unless otherwise stated, all spaces considered from now on are
infinite. Also, recall that \cite{engelking1989general} is our reference for topological cardinal functions.

In Corollary~\ref{cardinality_perfectly_normal_U_X} we were able to
calculate the precise value of $|\mathscr U_X|$ in terms of the
cardinal function $o(X)$, when $X$ belongs to
the class of $T_6$-spaces. Here, we present some
other classes of topological spaces in which the cardinality of the
lattice $\mathscr U_X$ can be determined in a similar fashion. 

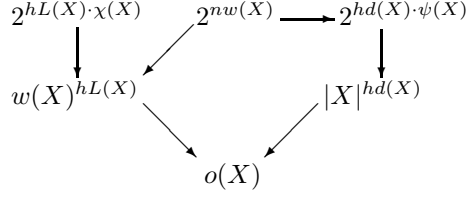
\begin{figure}[!ht]
  \begin{center}
    \setlength{\unitlength}{10pt}
    \begin{picture}(17,7)
      \put(0,6){$2^{hL(X)\cdot\chi(X)}$}
      \put(2.5,5.9){\vector(0,-1){1.9}}
      \put(0,3){$w(X)^{hL(X)}$}

      \put(5,3){\vector(1,-1){2}}
      \put(7.3,0){$o(X)$}
      \put(7,6){$2^{nw(X)}$}
      \put(6.9,6){\vector(-1,-1){1.9}}

      \put(10.2,6.2){\vector(1,0){2}}
      \put(12.4,6){$2^{hd(X)\cdot\psi(X)}$}

      \put(11.6,3){\vector(-1,-1){2}}
      \put(11.8,3){$|X|^{hd(X)}$}
      \put(14,5.8){\vector(0,-1){1.8}}
    \end{picture} 
  \end{center}
\caption{In this diagram $X$ is an arbitrary space and the symbol
$\kappa\to\lambda$ means that $\kappa\geq\lambda$. The upper right
inequality can be found in \cite[Theorem~7.1, p.~311]{hodel1978number}
and the rest of them are basic (see \cite{hodel1984cardinal}).}\label{topological_properties_bound_U_X}
\end{figure}

\begin{proposition}\label{cardinality_topological_properties_U_X}
  Given a space $X$, if any of the following statements holds, then $|\mathscr U_X|=2^{\mathfrak{c}}$.
  \begin{enumerate}
  \item $X$ is hereditarily Lindel\"of and first countable.
  \item $X$ admits a countable network.
  \item $X$ is hereditarily separable and has countable pseudocharacter. 
  \end{enumerate}
\end{proposition}
\begin{proof}
  From Proposition~\ref{roblemma} and
  Figure~\ref{topological_properties_bound_U_X} we deduce that $|\mathscr U_X|\leq 2^{\mathfrak{c}}$. The reverse inequality is a
  consequence of Corollary~\ref{cardinality_discrete_U_X}.
\end{proof}

In what follows, given a space $X$, we will employ the inequalities presented in
Figure~\ref{topological_properties_bound_U_X} together with
Proposition~\ref{roblemma}(2) in order to get bounds for $|\mathscr
U_X|$. 

Now, regarding compact spaces we have the following
results. 

\begin{lemma}\label{compact_bound}
  For any compact space $X$, $|\mathscr U_X| \leq \beth_2(hL(X))$.
\end{lemma}
\begin{proof}
  Given the hypotheses on $X$, we obtain $\chi(X)=\psi(X)\leq
  hL(X)$ and thus, the inequality needed follows from
  Figure~\ref{topological_properties_bound_U_X} and Proposition~\ref{roblemma}.
\end{proof}

\begin{proposition} 
  If $X$ is a compact space in which every open subset of it is an $F_\sigma$-set, then $|\mathscr U_X|=2^{\mathfrak{c}}$. In particular, every compact metrizable space satisfies the previous equality.
\end{proposition}
\begin{proof} It is sufficient to notice that our assumptions on $X$ imply $hL(X)=\omega$. Thus, Corollary~\ref{cardinality_discrete_U_X} and Lemma~\ref{compact_bound} give the desired result.
\end{proof}

Given an infinite cardinal $\kappa$, let us denote by $D(\kappa)$ and
$\beta D(\kappa)$ the discrete space of size $\kappa$ and its Stone-\v
Cech compactification, respectively. The regularity of $\beta
D(\kappa)$ implies that (see \cite[Theorem~3.3,
p.~11]{hodel1984cardinal}) \[w(\beta D(\kappa))\leq 2^{d(\beta
    D(\kappa))} = 2^\kappa.\] Therefore,
from Figure~\ref{topological_properties_bound_U_X} and the compactness
of $\beta D(\kappa)$ we deduce that \[|\mathscr U_{\beta D(\kappa)}|
  \leq \beth_2\left(nw\left(\beta D(\kappa)\right)\right) =
  \beth_2\left(w\left(\beta D(\kappa)\right)\right) \leq
  \beth_3(\kappa).\] On the other hand, since $\left|\beta
  D(\kappa)\right|=\beth_2(\kappa)$, Proposition~\ref{roblemma}(2)
gives \[\beth_3(\kappa)=2^{\left|\beta D(\kappa)\right|}\leq |\mathscr
  U_{\beta D(\kappa)}|.\] In conclusion, for any infinite cardinal
$\kappa$, $|\mathscr U_{\beta D(\kappa)}|=\beth_3(\kappa)$.

Once again, let $\kappa\geq\omega$ be a cardinal. If $D(2)$ is the
discrete space of size $2$, then $D(2)^\kappa$ is the Cantor cube of
weight $\kappa$. Clearly (see
Figure~\ref{topological_properties_bound_U_X}), \[|\mathscr
  U_{D(2)^{\kappa}}| \leq
  \beth_2\left(nw\left(D(2)^{\kappa}\right)\right) =
  \beth_2\left(w\left(D(2)^{\kappa}\right)\right) = \beth_2(\kappa).\]
Also, Proposition~\ref{roblemma}(2) produces \[\beth_2(\kappa)=2^{\left|D(2)^{\kappa}\right|}\leq |\mathscr
  U_{D(2)^{\kappa}}|.\] Hence, for any infinite cardinal
$\kappa$, $|\mathscr U_{D(2)^{\kappa}}| = \beth_2(\kappa)$.

Let $\mathbb{L}$ be the lexicographic square (i.e., $\mathbb{L}$ is
the cartesian product $[0,1]^{2}$ endowed with the topology generated
by the lexicographical ordering). By setting $Y :=
[0,1]\times\{\tfrac{1}{2}\}$ one gets a discrete subspace of $\mathbb
L$ and so, according to Corollaries~\ref{cardinality_perfectly_normal_U_X} and
\ref{order_embedding}, $\beth_2(\mathfrak c) = |\mathscr
U_Y|\leq|\mathscr U_{\mathbb L}|$. Finally, our definition of $\mathbb
L$ gives $o(\mathbb L)\leq 2^{\mathfrak c}$ and, as a consequence, $|\mathscr U_{\mathbb L}|\leq\beth_2(\mathfrak c)$. In other words, $|\mathscr U_{\mathbb L}| =
\beth_2(\mathfrak c)$.

The subspace $[0,1]\times\{0,1\}$ of $\mathbb L$ is called the double
arrow space and we will denote it by $\mathbb A$. Since the subspace
$(0,1)\times\{0\}$ of $\mathbb A$ is homeomorphic to Sorgenfrey's
line, the space $\mathbb A^2$ contains a discrete subspace of size $\mathfrak
c$. Therefore, as we did for $\mathbb L$, $|\mathscr U_{\mathbb
  A^2}|\geq\beth_2(\mathfrak c)$. For the reverse inequality note that
$o(\mathbb A^2)\leq o(\mathbb L^2)\leq 2^{\mathfrak c}$ and so,
$|\mathscr U_{\mathbb A^2}| = \beth_2(\mathfrak c)$.

A final note regarding $\mathbb A$ is pertinent. From (\ref{PicRio3_eqROyd}) and the
fact that $\mathbb A$ is separable, we deduce that $|\mathop{\rm
  RO}(\mathbb A^2)|\leq\mathfrak c$ and hence, \[\max\{2^{|\mathbb A^2|} ,
2^{|\mathop{\rm RO}(\mathbb A^2)|}\} = 2^{\mathfrak c} <
\beth_2(\mathfrak c) = |\mathscr U_{\mathbb A^2}|.\] This shows that
the lower bounds for $|\mathscr U_X|$ presented in
Proposition~\ref{roblemma}(2) need to be improved. 

\begin{proposition}\label{her_Lin} If $X$ is hereditarily Lindel\"of, then $|\mathscr U_X|=2^{o(X)}$.

\end{proposition}

\begin{proof} With Corollary~\ref{cardinality_perfectly_normal_U_X} in
  mind, we only need to show that all open subsets of $X$ are
  $F_\sigma$. Let $U\in \tau_X$ be arbitrary. For each $x\in U$ there
  is $U_x\in \tau_X$ such that $x\in U_x \subseteq \overline{U_x} \subseteq U$. Since
  $U$ is Lindel\"of, for some $F\in [U]^{\leq \omega}$ we obtain $U= \bigcup\left\{\overline{U_x} : x\in F\right\}$.
\end{proof}

We present now our findings regarding the following question.

\begin{question}\label{PicRio_cuandoU_XEs2AlaO}
  Given a space $X$, what conditions on $X$ imply that $|\mathscr U_X| = 2^{o(X)}$?
\end{question}

\begin{lemma}\label{PicRio_potenciaDehdX}
  If $X$ is a space with $|X|^{hd(X)} = |X|$, then $|\mathscr U_X|=2^{o(X)}$.
\end{lemma}
\begin{proof}
  It follows from Figure~\ref{topological_properties_bound_U_X} and
  our hypotheses that
  $o(X)\leq |X|$. On the other hand, the fact that $X$ is Tychonoff
  clearly implies the relation $|X|\leq o(X)$. Hence, the equality we
  need is a consequence of Proposition~\ref{roblemma}(2).
\end{proof}

\begin{proposition}\label{power_hd}
  If $X$ is a space for which there is a cardinal $\kappa$ with $|X|=2^{\kappa}$ and $\kappa \geq hd(X)$, then $|\mathscr U_X|=2^{o(X)}$.
\end{proposition}
\begin{proof}
  Our choice for $\kappa$ gives $|X|^{hd(X)} = |X|$ and so, the
  hypotheses of Lemma~\ref{PicRio_potenciaDehdX} are satisfied. 
\end{proof}

As usual, the acronym \textsf{GCH} stands for the Generalized
Continuum Hypothesis and $\mathop{\rm cf}(\alpha)$ denotes the
cofinality of an ordinal $\alpha$. 

\begin{proposition}\label{GCH_cofinality}
  Assuming {\sf GCH}, if $X$ is a space satisfying $\cf(|X|)>hd(X)$, then $|\mathscr U_X|=2^{o(X)}$.
\end{proposition}

\begin{proof}
  According to \cite[Lemma~10.42, p.~34]{kunen1980set}, $|X|^{hd(X)}=|X|$ and
  therefore we only need to invoke Lemma~\ref{PicRio_potenciaDehdX}.
\end{proof}

\begin{proposition}
  Given a space $X$, if $|X|$ is a singular strong limit cardinal, then $|\mathscr U_X|=2^{o(X)}$.

\end{proposition}

\begin{proof} The hypothesis allows us to use \cite[Theorem~3, p.~22]{hajnal1970discrete} to find a discrete set $D\subseteq X$ such that $|D|=|X|$. Hence, Proposition~\ref{roblemma}(2) and Corollary~\ref{cardinality_discrete_U_X} imply that $|\mathscr U_X|=2^{o(X)}$.
\end{proof}

Let us denote by {\sf A} the statement \lq\lq {\sf GCH} holds and
there are no inaccessible cardinals.{\rq\rq}

\begin{corollary}
  Assume {\sf A} holds. Then, for any space $X$ whose cardinality is a
  limit cardinal we obtain $|\mathscr U_X|=2^{o(X)}$.

\end{corollary}

With the idea in mind of finding the effect that {\sf GCH} has on
$|\mathscr U_X|$, let us recall that, for a cardinal number $\kappa$, $\kappa^{+}$ represents the successor cardinal of $\kappa$.

\begin{proposition}
  If {\sf GCH} holds, then, for any space $X$, $|\mathscr U_X|$ is a
  regular uncountable cardinal. 
\end{proposition}

\begin{proof} On the one hand,
  Corollary~\ref{cardinality_discrete_U_X} implies that
  $|\mathscr{U}_X|$ is uncountable. On the other hand, since $2^{|X|}
  \leq |\mathscr U_X| \leq 2^{o(X)} \leq \beth_2\left(|X|\right) =
  (2^{|X|})^+ $, we deduce that $|\mathscr U_X| \in \{|X|^{+},
  (2^{|X|})^+ \}$. In either case, $|\mathscr U_X|$ is regular. 
\end{proof}

\begin{proposition}\label{GCH_her_separable} Under the assumptions
  $\mathfrak c = \omega_1$ and $2^{\mathfrak c}=\omega_2$, if $X$ is a
  hereditarily separable space, then $|\mathscr U_X|=2^{o(X)}$.

\end{proposition}

\begin{proof}
  According to \cite[Theorem~4.12, p.~21]{hodel1984cardinal}, the
  relation $hd(X)=\omega$ guarantees that $|X| \leq 2^{\mathfrak{c}}$
  and consequently, $|X| \in \{\omega, \omega_1, \omega_2\}$.

  When $|X| \in \{\omega_1, \omega_2\}$, Proposition~\ref{power_hd}
  gives us the desired equality. Finally, if $|X|=\omega$, then $X$
  admits a countable network and thus (see
  Proposition~\ref{cardinality_topological_properties_U_X}),
  $|\mathscr U_X| = 2^{\mathfrak c} = 2^{o(X)}$.
\end{proof}

Suppose $X$ is a space. Since $\mathscr U_X$ is a subset of
$\Sigma(C(X))$, we obtain $|\mathscr U_X| \leq |\Sigma(C(X))|$. With
the idea in mind of showing two examples for which this inequality is strict, let us
note first that the fact $|C(X)|\geq\omega$ implies, according to \cite[Theorem~1.4,
p.~179]{larson1975lattice}, that $|\Sigma(C(X))|=\beth_2(|C(X)|)$.

When $X$ is an infinite discrete space, we obtain $|C(X)| = 2^{|X|}$
and so, by Proposition~\ref{roblemma}(2), \[|\mathscr U_X| \leq
  \beth_2(|X|)<\beth_3(|X|) = \beth_2(|C(X)|).\]

On the other hand, if $X$ is any infinite countable space, then it follows from Proposition~\ref{cardinality_topological_properties_U_X}(2) that \[|\mathscr U_X| = 2^{\mathfrak{c}} < \beth_2(\mathfrak{c}) \leq \beth_2(|C(X)|).\]

Our final result of this section establishes some conditions for a
family of topological spaces under which the corresponding Tychonoff
product $X$ satisfies the equality $|\mathscr U_X| =
|\Sigma(C(X))|$. For this proposition we won't require for our spaces
to be infinite.

\begin{proposition}
  Assume that $\kappa$ is an infinite cardinal. Let $X$ be the topological
  product of a family of spaces $\{X_\xi : \xi<2^{\kappa}\}$. If $|X_\xi|\geq2$ and $d(X_\xi)\leq
  \kappa$ for each $\xi<\kappa$, then $|\mathscr U_X| = |\Sigma(C(X))|$.
\end{proposition}

\begin{proof} Since we always have the inequality $|\mathscr U_X| \leq
  |\Sigma(C(X))|$, we only need to show that $|\mathscr
  U_X|\geq\beth_2(|C(X)|)$.

  According to Proposition~\ref{roblemma}(2), $|\mathscr U_X|\geq
  2^{|X|}$. Now, the fact that each $X_\xi$ has at least two points
  gives $|X|\geq \beth_2(\kappa)$ and so, $2^{|X|}\geq
  \beth_3(\kappa)$. On the other hand, the Hewitt-Marczewski-Pondiczery Theorem (see
  \cite[Theorem~2.3.15, p.~81]{engelking1989general}) implies that $d(X)\leq
  \kappa$ and therefore, from the well-known relation
  $2^{d(X)}\geq|C(X)|$ we deduce that $2^\kappa\geq|C(X)|$. In
  conclusion, $|\mathscr U_X|\geq\beth_3(\kappa)\geq\beth_2(|C(X)|)$, as required. 
\end{proof}

For example, if $X$ is a Cantor cube of the form $D(2)^{2^{\kappa}}$, where
$\kappa$ is an infinite cardinal, then $|\mathscr
U_X|\geq\beth_2(|C(X)|)$.

We close the paper with a list of open questions.

\begin{question}
  Does Corollary~\ref{cardinality_perfectly_normal_U_X} remain true if
  we replace $T_6$ with $T_5$ in the hypotheses?
\end{question}

\begin{question}
  Regarding Proposition~\ref{her_Lin}, is it true that for any compact
  space $X$, $|\mathscr U_X|= 2^{o(X)}$?
\end{question}

\begin{question}
  Can we drop the set-theoretic assumptions $\mathfrak c = \omega_1$
  and $2^{\mathfrak c}=\omega_2$ in Proposition~\ref{GCH_her_separable}?
\end{question}

We conjecture that, under \textsf{A}, the equality
\begin{align}
  \label{PicRioeq:cardDeoX}
  |\mathscr U_X| = 2^{o(X)}
\end{align}
 holds for any space $X$. Even though we did not prove or
refute this conjecture, we were able to obtain some partial results
(for example, if one assumes \textsf{A}, then (i) for
any space $X$, $\beth_2(s(X))\leq|\mathscr U_X|$, and (ii) we possess a
short list of classes $\mathscr S$ in such a way that $X\in\mathscr S$
implies that (\ref{PicRioeq:cardDeoX}) holds). Consequently, we pose the following
problem. 

\begin{question} Does it follow from \textsf{A} that
  (\ref{PicRioeq:cardDeoX}) is true for any space $X$?
\end{question}


\begin{thebibliography}{5}
  \thispagestyle{myheadings}
  
   \bibitem{engelking1989general} R. Engelking, \textit{General Topology}, second ed., Sigma Series in Pure Mathematics, vol. 6, Heldermann Verlag, Berlin, 1989.
   
   \bibitem{hajnal1970discrete} A. Hajnal, I. Juh\'asz, \textit{Discrete subspaces of topological spaces, \textsc{II}}, Indag. Math., {\bf 71(1)} (1970), 18--30.  
  
   \bibitem{hodel1984cardinal} R. Hodel, \textit{Cardinal Functions I}, Handbook of Set-Theoretic Topology (K. Kunen and J. E. Vaughan, eds.), Amsterdam, pp. 1--61, 1984.
  
   \bibitem{hodel1978number} {\bf \textemdash}, \textit{The number of closed subsets of a topological space}, Canadian Journal of Mathematics, {\bf 30}(2) (1978), 301--314. 

   \bibitem{jech2003set} T. Jech, \textit{Set Theory. The third millenium edition, revised and expanded}, Springer Monograph in Mathematics, Springer-Verlag Berlin Heidelberg, 2003.

   \bibitem{koppelberg1989handbook1} S. Koppelberg, \textit{General Theory of Boolean Algebras}, in \textit{Handbook of Boolean algebras} (eds. J.D. Monk and R. Bonnet), North-Holland, Amsterdam, 1989.

   \bibitem{kunen1980set} K. Kunen, \textit{Set theory. An Introduction to Independence Proofs}, Studies in Logic and the Foundations of Mathematics, vol. 102, North-Holland Publishing Co., Amsterdam, 1980. 

   \bibitem{larson1975lattice} R.E. Larson and J.A. Susan, \textit{The lattice of topologies: A survey}, The Rocky Mountain Journal of Mathematics, {\bf 5}(2) (1975), 177--198.

   \bibitem{perez2021lattice} L.A. P\'erez-Morales, G. Delgadillo-Pi\~n\'on, and R. Pichardo-Mendoza, \textit{The lattice of uniform topologies on $C(X)$}, Questions and Answers in General Topology {\bf 39} (2021), 65--71.

   \bibitem{pichardo2013pseudouniform} R. Pichardo-Mendoza, \'A. Tamariz-Mascar\'ua, and H. Villegas-Rodr\'{\i}guez, \textit{Pseudouniform topologies on $C(X)$ given by ideals}, Comment. Math. Univ. Carolin. {\bf 54}(4) (2013), 557--577.

\end{thebibliography}
\end{document}